\documentclass[a4paper, BCOR=0.0mm, DIV=calc]{scrartcl}
\usepackage[T1]{fontenc}
\usepackage[utf8]{inputenc}
\usepackage{url}
\usepackage{lmodern}
\usepackage{microtype}
\usepackage[english]{babel}
\usepackage{amsmath, amssymb, amsfonts}
\usepackage{nicefrac}
\newcommand{\diff}[1]{\mathrm{#1}}
\newcommand{\gam}[1]{\Gamma ( #1 )}

\KOMAoptions{twoside=false, twocolumn=false, headinclude=false, footinclude=false, mpinclude=false, pagesize=auto}
\recalctypearea

\begin{document}
\title{Note on Malmstèn's paper De Integralibus quibusdam definitis seriebusque infinitis }
\author{Alexander Aycock}
\date{}
\maketitle
\begin{abstract}
We present a proof of the functional equation of the Riemann zeta-function or more precisely the Dirichlet eta-function, which proof seems to be new but follows almost immediately from Malmstèn's paper ``De integralibus quibusdam definitis seriebusque infinitis''\cite{M}.
\end{abstract}
\section*{Introduction}
In 1846 the Swedish mathematician C. J. Malmstèn wrote the paper ``De integralibus quibusdam definitis seriebusque infinitis'' \cite{M}, in which he demonstrates some very remarkable formulas, like
\[
\int_0^1 \frac{\ln{\ln{\frac{1}{x}}}}{1+x^2}\diff{d}x  = \frac{\pi}{2}\log{\frac{\gam{\frac{3}{4}}}{\gam{\frac{1}{4}}}}\sqrt{2\pi}
\]
It seems that the paper was completly forgotten, as the above formula
was rediscovered $140$ years later, \cite{V}. The book {\em Handbuch der Theorie der Gammafunktion} \cite{N} contains a reference to it.
The highlight of the paper is certainly the derivation of the
functional equation for some Dirichlet series. He shows for example
\[
L(1-s) = \left(\frac{2}{\pi}\right)^s\sin{\left(\frac{\pi s}{2}\right)}\gam{s}L(s)  \quad \text{with} \quad L(s) = \sum_{n=0}^{\infty}\frac{(-1)^n}{(2n+1)^s}
\]
which is often called Dirichlet's beta-function after Dirichlet and denoted by $\beta(s)$ nowadays. In the paper one finds several other functional equations, only the one of the Riemann zeta-function 
\[\zeta{(s)} = \sum_{n=1}^{\infty}\frac{1}{n^s}\]
or equivalently the Dirichlet eta-function 
\[\eta{(s)} = \sum_{n=1}^{\infty}\frac{(-1)^{n+1}}{n^s}\]
is missing. Reading through this rich paper attentively, one discovers that Malmstèn could also have shown these functional equations with
equal ease. The purpose of this little note is show this and as such it stands in the good tradition of giving "missed proofs", like \cite{P}.
We will focus on giving a proof of the formula
\[
\eta(1-s) = \frac{2^s-1}{1-2^{s-1}}\pi^{-s}\cos{\left(\frac{\pi s}{2}\right)}\gam{s}\eta(s)
\]
Our starting point is formula 30 in Malmstèn's paper
\[
\int_0^\infty\frac{e^{au}-e^{-au}}{e^{\pi u}-e^{-\pi u}}\cdot \frac{\diff{d}u}{u^s}=\frac{\sin{(a)}}{\gam{s}\cos{\frac{\pi s}{2}}}\int_0^1\frac{\ln^{s-1}{(\frac{1}{y})}}{1+2y\cos{a} + y^2}\diff{d}y\;.
\]
This formula is quite easy to prove by expanding the function on the left hand side into an infinite series and integrating term-by-term, what will lead to a well-known partial fraction decomposition and then almost immediately to the right hand side of the equation. \\

We now divide both sides by $\sin{a}$ and use the limit 
\[
\lim_{a \rightarrow 0}\frac{e^{au}-e^{-au}}{\sin{a}} = 2u
\]
to obtain
\[
2\cos{\left(\frac{\pi s}{2}\right)}\int_0^{\infty}\frac{u^{1-s}}{e^{\pi u}-e^{-\pi u}}\diff{d}u = \frac{1}{\gam{s}}\int_0^1\frac{\ln^{s-1}{(\frac{1}{y})}}{(1+y)^2}\diff{d}y
\]
and therefore
\[
2\cos{\left(\frac{\pi s}{2}\right)}\int_0^{\infty}\frac{u^{1-s}e^{-\pi u}}{1-e^{-2\pi u}}\diff{d}u = \frac{1}{\gam{s}}\int_0^1\frac{\ln^{s-1}{(\frac{1}{y})}}{(1+y)^2}\diff{d}y
\]   
The first formula in fact can also be found in Malmstèn's paper as formula 57 for $x = 0$. And now I claim, that in a sense this is already the functional equation for $\eta{(s)}$!

We now start with the proof. Using the series expansions for $\frac{1}{1-e^{-2\pi u}}$ and $\frac{1}{(1+y)^2}$ we obtain
\[
2\cos{\left(\frac{\pi s}{2}\right)}\int_0^{\infty}\sum_{n=0}^{\infty}e^{-(2n+1)\pi u}u^{1-s}\diff{d}u = \frac{1}{\gam{s}}\int_0^1\ln^{s-1}{\left(\frac{1}{y}\right)\sum_{n=1}^{\infty}(-1)^{n+1}n y^{n-1}\diff{d}y}
\]
By interchanging the order of summation and integration this becomes
\begin{align*}
&2\cos{\left(\frac{\pi s}{2}\right)}\sum_{n=0}^{\infty}\int_0^\infty e^{-(2n+1)\pi u}u^{1-s}\diff{d}u \\
=& \frac{1}{\gam{s}}\sum_{n=1}^{\infty}(-1)^{n+1}n\int_0^1\ln^{s-1}{\left(\frac{1}{y}\right)y^{n-1}\diff{d}y}
\end{align*}
Writing  $(2n+1)\pi u = y$ in the $n$-term of the left hand side, we obtain
\begin{align*}
&2\cos{\left(\frac{\pi s}{2}\right)}\sum_{n=0}^{\infty}\int_0^{\infty}e^{-y}\left(\frac{y}{(2n+1)\pi}\right)^{1-s}\frac{\diff{d}y}{(2n+1)\pi} \\
=& \frac{1}{\gam{s}}\sum_{n=1}^{\infty}(-1)^{n+1}n\int_0^1\ln^{s-1}{\left(\frac{1}{y}\right)y^{n-1}\diff{d}y}
\end{align*}
Both integrals are standard and can be expressed in terms of the $\Gamma$-function:
\[\int_0^{\infty}e^{-y}y^{1-s}\diff{d}y={\gam{2-s}}\]
and
\[\int_0^1\ln^{s-1}\left({\frac{1}{y}}\right) y^{n-1}\diff{d}y = \frac{\gam{s}}{n^s}\]
so we obtain
\[
\frac{2\cos{\left(\frac{\pi s}{2}\right)}}{\pi^{2-s}}\sum_{n=0}^{\infty}\frac{1}{(2n+1)^{2-s}}{\gam{2-s}} 
= \frac{1}{\gam{s}}\sum_{n=1}^{\infty}(-1)^{n+1}n\frac{\gam{s}}{n^s}\;.
\]
Using the notation
\[
 \lambda{(s)} = \sum_{n=0}^{\infty}\frac{1}{(2n+1)^s}
\]
this can be written as
\[
\frac{2\cos{\left(\frac{\pi s}{2}\right)}}{\pi^{2-s}}\lambda{(2-s)}\gam{2-s} = \eta{(s-1)}
\]
But note the relation between $\lambda(s)$ and $\eta(s)$ which is elementary to 
prove:
\[ \lambda(s) = \frac{2^s-1}{2^s-2} \eta(s)\]
We now use this to replace $\lambda{(2-s)}$ with the  corresponding  $\eta{(2-s)}$ to get
\[
\frac{2\cos\left(\frac{\pi s}{2}\right)}{\pi^{2-s}} \frac{2^{2-s}-1}{2^{2-s}-2}\eta{(2-s)}\gam{2-s} = \eta{(s-1)}
\]
Canceling the 2 and writing $s$ instead of $2-s$ we find
\[
\eta{(1-s)} = \frac{2^s-1}{1-2^{s-1}}\pi^{-s}\cos{\left(\frac{\pi s}{2}\right)}\gam{s}\eta{(s)}
\]
the functional equation for $\eta$!\\

As one sees from this derivation,  we have used only formulas well known to Malmstèn and one wonders why Malmstèn missed 
this formula. Of course, a precise discussion of the domains of validity of the given formulas and the idea of analytic continuation are not to be found explicitly in Malmstèn paper, but nevertheless the would have led him to the right equation. We note that Riemann's celebrated paper \cite{R} that contains the first rigorous proof of the functional equation for $\zeta$ is from $1859$. 

But Malmstèn uses formula 30 to obtain the Fourier series for $\ln\gam{(x)}$ which is not related to $\eta(s)$ directly and
from that vista the functional equation is no longer in sight. He finds the Fourier series of $\ln\gam{x}$ in a very clever 
way without calculating the coefficients explicitly and represents it as follows
\[
	\int_0^1 \frac{ \ln{\ln{(\frac{1}{y})}} \diff{d}y }{ 1 + 2y\cos{a} + y^2} = \frac{\pi}{2\sin{a}}\ln\left\{ \frac{(2\pi)^{\frac{a}{\pi}} \Gamma\left( \frac{1}{2} + \frac{a}{2\pi}\right)}{\Gamma \left( \frac{1}{2} - \frac{a}{2\pi} \right)} \right\}
\]
which is easily seen to be equivalent to the expression Kummer \cite{K} gave one year later in 1847.\\

There is a lot more to say about this very interesting paper that seems not very well known. Although it appeared in
{\em Journal f\"ur Mathematik} it was apparently not quoted by contemporary mathematicians. Apparently, the paper was
forgotten, one of the reasons being maybe that it was only available in Latin.  But I think that the functional equation of 
the zeta-function, which is equivalent to the one we gave here, illustrates very well that it contains some very intriguing and deep results. And Malmstèn obtains them in an quite unusual way. 
We intend to review the contents of this paper in more detail on another occasion, on which we will also provide a translation 
of the Latin original.\\[2mm]

{\bf Acknowledgement:} I would like to thank Gert Almkvist and Duco van Straten for showing interest in this note.

\end{document}